\documentclass[a4paper,10pt]{article}
\usepackage[utf8x]{inputenc}
\usepackage[russian]{babel}

\usepackage{amssymb}

\title{Towards the moduli space of special Bohr - Sommerfeld lagrangian cycles }
\author{{\bf Nikolay A. Tyurin}\footnote{BLTPh JINR (Dubna) and LAG NRU HSE (Moscow),  {\bf ntyurin@theor.jinr.ru}.
The study has been funded by the Russian Academic Excellence Project '5-100'.}}

\begin{document}

\maketitle

\begin{abstract} In previous papers we introduced the notion of special Bohr - Sommerfeld lagrangian cycles on a compact simply connected symplectic manifold
with integer symplectic form, and presented the main interesting case: compact simply connected algebraic variety with an ample line bundle such that
the space of Bohr - Sommerfeld lagrangian cycles with respect to a compatible Kahler form of the Hodge type and holomorphic sections of the bundle is finite.
The main problem appeared in this way is singular components of the corresponding lagrangian shadows (or sceletons of the corresponding Weinstein domains)
which are hard to distinguish or resolve. In the present text we avoid this difficulty
presenting the points of the moduli space of special Bohr - Sommerfeld lagrangian cycles by exact compact lagrangian submanifolds on the complements $X \backslash D_{\alpha}$
modulo Hamiltonian isotopies, where $D_{\alpha}$ is the zero divisor of holomorphic section $\alpha$. In a sense it corresponds to the usage of gauge classes of hermitian connections instead of
pure holomorphic structures in the theory of the moduli space of (semi) stable vector bundles.
\end{abstract}

\section{General theory}

Consider $(M, \omega)$ --- a compact simply connected symplectic manifold of dimension $2n$, endowed with a symplectic form of integer type, $[\omega] \in H^2(M, \mathbb{Z})$.
Then it exists a prequantization data --- the pair $(L, a)$, where $L \to M$ is a hermitian line bundle and $a \in {\cal A}_h(L)$ is a hermitian connection such that the curvature form $F_a = 2 \pi i \omega$
(thus the first Chern class $c_1(L) = [\omega]$). 

An $n$ - dimensional submanifold $S \subset M$ is called {\it lagrangian} iff the restriction $\omega|_S$ identically vanishes; $S$ is called {\it Bohr - Sommerfeld} lagrangian (or BS for short)
iff the resctriction $(L, a)|_S$ admits a covariantly constant section $\sigma_S \in \Gamma (L|_S)$, defined up to $\mathbb{C}^*$. For any choosen smooth section
$\alpha \in \Gamma(M, L)$ we say that $S \subset M$ is special with respect to $\alpha$ Bohr - Sommerfeld lagrangian cycles (or $\alpha$ - SBS for short) iff it is 
Bohr - Sommerfeld lagrangian and the restriction
$\alpha|_S = e^{i c} f \sigma_S$, where $c$ is a real constant and $f $ is a strictly positive real function on $S$. In the present paper we consider compact orientable lagrangian submanifolds
only.

It was already shown that the definition above can be reformulated in terms of calibrated lagrangian geometry. For any smooth section $\alpha \in \Gamma (M, L)$ we define
the complex valued 1 -form
$$
\rho_{\alpha} = \frac{<\nabla_a \alpha, \alpha>}{<\alpha, \alpha>} \in \Omega^1_{\mathbb{C}}(M \backslash D_{\alpha})
$$
where $D_{\alpha} = \{ \alpha = 0 \} \subset M$ is the zeroset of $\alpha$. This form satisfies the floowing properties: its real part is exact being $d ({\rm ln} \vert \alpha \vert)$,
and the imaginary part is a canonical 1- form on the complement $M \backslash D_{\alpha}$ since $d ({\rm Im} \rho_{\alpha}) = 2 \pi \omega$.

In these terms an $n$ - dimensional submanifold $S \subset M$ is  $\alpha$ - SBS lagrangian if and only if the restriction ${\rm Im} \rho_{\alpha}|_S$
identically vanishes (the proof and details can be found in [1]).

Using this ``calibrated reformulation'' of the definition  one proved that any Weinstein neighborhood ${\cal O}(S_0)$ of an  $\alpha$- SBS lagrangian submanifold $S_0$ cann't contain any other
$\alpha$ - SBS lagrangian submanifold of the same type. It follows that a fixed $\alpha$ admits a discrete set of $\alpha$ - SBS lagrangian submanifolds of the same topological type.

Recall that  the situation stated above is the input of ALAG - programme, proposed by A. Tyurin and A. Gorodentsev in [2]: starting with such $(M, \omega)$ they constructed
certain moduli space of Bohr - Sommerfeld lagrangian cycles of fixed topological type, denoted as ${\cal B}_S$. Such a moduli space is a Frechet smooth infinite dimensional
real manifold, locally modelled by unobstracted isodrastic deformations  of BS lagrangian submanifolds. To define ${\cal B}_S = {\cal B}_S({\rm top} S, [S])$ one has 
to choose topological type of $S$ and the homology class
$[S] \in H_n(M, \mathbb{Z}$ of the corresponding BS submanifolds. Moreover, the BS - level can be shifted up, so one has a series of the moduli space ${\cal B}_S^k$
(details see in [2]).
 
Therefore in the situation presented above  we can consider in the direct product ${\cal B}_S \times \mathbb{P} \Gamma (M, L)$ certain subset given by the condition: pair $(S, p) \in {\cal U}_{SBS}$
iff $S$ is $\alpha$ - SBS lagrangian submanifold where $\alpha$ corresponds to point $p$ in the projectivized space (and of course it is possible to shift the BS - level, getting
the corresponding subset in the direct product ${\cal B}_S^k \times \mathbb{P} \Gamma (M, L^k)$, but in the present text we leave aside the variation of BS- level). 

This subset ${\cal U}_{SBS}$ was studied in [1]; the main result is that the canonical projection $p: {\cal U}_{SBS} \to \mathbb{P} \Gamma (M, L)$ has discrete fibers,
has non degenerated differential in smooth points and projects ${\cal U}_{SBS}$ to an open subset of the last projective space. As a corollary one establishes that
${\cal U}_{SBS}$ admits a Kahler structure at smooth points. It seems to be interesting since we have started from pure symplectic situation and came to
an object from the Kahler geometry.

\section{The case of algebraic varieties}

Let  $X$ be a compact smooth simply connected  algebraic variety which admits an ample line bundle $L$; then it can be regarded as a special case of the situation presented above.

Indeed, fixing an appropriate hermitian structure $h$ on $L$ one induces the corresponding Kahler form $\omega$: any holomorphic section $\alpha \in H^0(X, L)$
in the presence of $h$ defines the function $\psi_{\alpha} = - {\rm ln} \vert \alpha \vert_h$ on the complement $X \backslash D_{\alpha}$ which is a Kahler potential, therefore
$\omega$ is given by $ d I d \psi_{\alpha}$, and the ampleness condition ensures that whole $X$ is covered by the complements to divisors from the complete linear system $\vert L
\vert = \mathbb{P} H^0(X, L)$, so $\omega$ is globally defined in $X$, see [3].

Thus one can consider  $(X, L)$ as a symplectic manifold with integer symplectic form, and $L$, the prequantization line bundle,  is automatically endowed with a
prequantization connection $a$, compatible with the holomorphic structure on the bundle.
For a holomorphic section $\alpha$ one has $\nabla_a \alpha = \partial_a \alpha$ and consequently the form $\rho_{\alpha}$ has type (1, 0) with respect to the complex structure.
Then one can deduce that  the SBS condition with respect to a holomorphic section is equivalent to the following condition: a lagrangian submanifold  $S \subset X$ is 
$\alpha$- SBS if and only if it is invariant under
the flow generated by the gradient vector field ${\rm grad} \psi_{\alpha}$ (see [4]). 

It is well known is algebraic geometry fact: the complement $X \backslash D_{\alpha}$, described above, is an example of the Stein variety,
and since we would like to study lagrangian geometry of these complements we must follow the key points of the programme ``From Stein to Weinstein and back'', see [5].
The situation we are studying here must be regarded in the framework of the Weinstein manifolds and Weinstien structures, see [5] and [6].
Indeed, the gradient vector field ${\rm grad} \psi_{\alpha}$ is  Liouville, while the function $\psi_{\alpha}$ is the second ingredient of the Weinstein structure
(of course, it just reflects the fact that $X \backslash D_{\alpha}$ is Stein). 

Since we claim that a lagrangian $S \subset X \backslash D_{\alpha}$ is $\alpha$ - SBS if and only if it is stable with respect to the gradient flow of
$\psi_{\alpha}$ it follows that such an $S$ must be contained by  the base set $B_{\alpha} \subset X \backslash D_{\alpha}$ defined as the union of (1) finite critical points of $\psi_{\alpha}$
and (2) finite trajectories of the gradient flow. Now we can translate our $\alpha$ - SBS condition to the language of Weinstein manifolds and structures:
a lagrangian submanifold $S \subset X$ is $\alpha$- SBS iff it is a component of the lagrangian sceleton defined by the Weinstein structure given by $({\rm grad} \psi_{\alpha}, \psi_{\alpha})$
on the complement   $X \backslash D_{\alpha}$.

{\bf Remark.} In the previous texts [4]  we use the term ``Lagrangian shadow of ample divisor'' for the lagrangian components of the lagrangian sceleton, 
since we would like to emphersize  the fact that the corresponding lagrangian components
arises for any ample divisor; in the theory of Weinstein manifolds which covers much wider situation than our complements $X \backslash D_{\alpha}$ one speaks about
regular lagrangian submanifolds. Below we use this parallel for the modified definition of moduli space of special Bohr - Sommerfled lagrangian cycles.

The old definition (see [4])  we have tried to exploite for the consruction of  certain moduli space of SBS lagrangian cycles over algebraic varieties used to be the following.
 Take the canonical projection $p: {cal U}_{SBS} \to \mathbb{P} \Gamma (M, L)$ to the second direct summand from the section 1, then in the present situation it is a finite dimensional
 projective subspace $\mathbb{P} H^0(X, L) \subset \mathbb{P} \Gamma(X, L)$ which corresponds to holomorphic sections.
Then the preimage ${\cal M}_{SBS} = p^{-1}(\mathbb{P} H^0(X, L))$ must be finite (and we have proved it for smooth Bohr - Sommerfeld
submanifolds in [4]), and we would like to understand it as the ``moduli space''. 

But the great problem appears in this case since  the componnents of the lagrangian sceleton (or the base set) $B_{\alpha}$ are very far from being smooth lagrangian submanifolds (or even cycles),
 therefore strictly speaking our coarse ``moduli space'' must be empty in many cases, and the framework of algebraic geometry doesn't admit any variational freedom to resolve this trouble. 
In the simple case, when $H_n(X \backslash D_{\alpha}, \mathbb{Z}) = \mathbb{Z}$ for generic smooth element $D_{\alpha}$ of the complete linear system $\vert L \vert$ the moduli space ${\cal M}_{SBS}$
 can be however correctly defined,
as it was done in [4], but in more geometrically interesting  cases we face great problem in this way: we must either present a strong theory of desingularization of
the components of lagrangian shadows doing it however in concordance with the technical details of ALAG or find a different definition
of special Bohr - Sommerfeld cycles with respect to holomorphic sections such that these new special submanifolds should be automatically smooth.

\section{Desingularizing the definition}

Recall that we study the lagrangian geometry of the complement $X \backslash D$, where $D$ is a compact smooth simply connected algebraic variety and $D$ is an ample divisor;
then we have fixed an appropriate hermitian structure $h$ on the ample line bundle $L \to X$, corresponding to $D$, and get the Kahler form $\omega$, such that
the function $\psi_{\alpha} = - {\rm ln} \vert \alpha \vert$ is a Kahler potential ($D = \{ \alpha = 0 \} \subset X$).

The Kahler potential $\psi_{\alpha}$ defines the structure of the Weinstein domain on $X \backslash D$, given by 1- form $\lambda = I d \psi_{\alpha}$ and $\psi_{\alpha}$ itself (see [5]);
then we can study {\it exact} compact orientable Lagrangian submanifolds in $X \backslash D_{\alpha}$, so is such  Lagrangian submanifolds $S \subset X \backslash D$ that the resctriction
$\lambda|_S$ is an exact form. Remark that any such an exact $S$ must be Bohr - Sommerfeld in whole $X$ with respect to the corresponidng prequantization data. Moreover,
we can introduce certain condition on Bohr - Sommerfeld lagrangian submanifolds in $X$ which is equivalent to the exactness condition on the complement $X \backslash D_{\alpha}$.

Namely, let $X \supset D$ is as above, and the corresponding symplectic form $\omega$   evidently rerpesents the cohomology class Poincare dual to $[D] \in H_{2n-2}(X, 
\mathbb{Z})$. Then we say that a lagrangian submanifold $S \subset X$ is {\it D - monotonic} iff $D \cap S = 0$ and for any oriented  loop $\gamma \subset S$ and any
compatible oriented disc $K_{\gamma} \subset X$, bounded by $\gamma$, the topological sum of the intersection points $D \cap K_{\gamma}$ equals to
the symplectic area of $K_{\gamma}$ (note that if $K_{\gamma}$ intersects $D$ non transversally then we can deform it to have transversal intersection). Of course, this
definition is applicable in much wider situation than stated above. 

Now it is not hard to see that 

{\bf Proposition.} {\it A lagrangian submanifold $S \subset X \backslash D$ is exact with respect to $\lambda$ if and only if $S$ is D - monotonic with respect to
$D$.}

In [6] one presents the list of open problems stated in the theory of Weinstein manifolds; and one of these problems hints how the definition
of spceial Bohr - Sommerfeld lagrangian submanifolds can be modified. Namely the Problem 5.1 from [6] asks are there non - regular exact lagrangian submanifolds
in $X \backslash D$? Regularity here means that they appear as components of lagrangian sceleton of the Weinstein domain; at the same time in our language regularity means that
they appear as components of tha Lagrangian shadow $Sh^{Lag}(D)$. If the answer is negative then we should get  desingularizations of the components
of $Sh^{Lag}(D)$ (or lagrangian sceleton) given by exact lagrangian submanifolds in the complement $X \backslash D$. Of course the space of exact lagrangian
submanifolds is too huge for our purposes (finitness of the moduli space), but we can factorize the space of all exact lagrangian submanifolds modulo
Hamiltonian isotopies.

Even  if Problem 5.1 has positive solutions we introduce the following

{\bf Definition.} {\it In the situation presented above consider the space of pairs $\tilde {\cal M}_{SBS} = \{([S], D)| S \in {\cal B_S}, D \in \vert L \vert \}$
where $[S]$ is a class of smooth compact orientable exact lagrangian submanifolds on the complement $X \backslash D$ up to Hamiltonian isotopy.}

The space $\tilde{\cal M}_{SBS}$ admits the forgetfull map $p: \tilde{\cal M}_{SBS} \to \vert L \vert$, and we can prove that the fibers are discrete
and that the differential of this map is non degenerated (the arguments are essentially the same as in [1]: we study the local picture over
a Weinstein neighborhood of a fixed D - monotonic $S \subset M$). 

In this setup the negative answer to Problem 5.1 from [6] should mean that we have exactly the same spaces: ${\cal M}_{SBS} = \tilde{\cal M}_{SBS}$
when the first space can be correctly defined (f.e. if $H_n(X\backslash D, \mathbb{Z}) = \mathbb{Z})$). At the same time
even if there are examples of non regular exact lagrangian submanifolds we still have some freedom to claim that the identity takes place.
For example, let it be a lagrangian sphere which appears in the case when one doesn't have a sphere as a regular component in the sceleton
(the mostly possible example for Problem 5.1), --- but since we've fixed the topological type ${\rm top} S$ to construct the moduli space ${\cal B}_S$
then take
it different from the spherical one, and then the choice ``kills'' inapproriate components. On the other hand if Problem 5.1 has negative answer
then it ensures that the forgetfull map $p$ has finite fibers, so the moduli space admits the structure of a finite covering of an open part
of the projective space.  

Let us illustrate  the story by the example which has appeared several times in the previous texts, see [arX]. Take $X = \mathbb{C} \mathbb{P}^1$ and
consider $L = {\cal O}(3)$. Study  the situation for certain concrite holomorphic section  f.e. for the section defined  by the polynomial $P_3 = z_0^3 -z_1^3$.
It vanishes at three roots of unity $p_1, p_2, p_3$ which become poles for the function $\psi = - {\rm ln} \vert P_3 \vert$; the last one has
exactly 5 finite critical points --- 2 local minima $m_1 =  [1:0], m_2 = [0:1],$ and three saddle points $s_1, s_2, s_3$ at the roots of -1. The base set
consists of three  lines $\gamma_i$ each of which joins $m_1$ and $m_2$ passing through $s_i$. Totally we get non smooth simple loops only in the base set:
each closed loop  is formed by two lines $\gamma_i, \gamma_j$, and at the points $m_1, m_2$ the loop has corners. Therefore if we are looking for the ``old version''
of the moduli space ${\cal M}_{SBS}$ we must specialize what singular loops are allowed in our situation. However in this case the specialization can be done:
we may say that a singular loop is allowed if it can be transformed by a small deformation to a smooth Bohr - Sommerfeld loop. Then one gets exactly three
simple elements for the moduli space. 

But our new ``desingularized'' definition of the moduli space works much better:
we  claim that there are exactly three smooth exact closed loops on the complement $\mathbb{C} \mathbb{P}^1 \backslash \{p_1, p_2, p_3\}$ up to Hamiltonian isotopy.
Indeed, for each zero $p_i$ we can take a smooth loop surrounding $p_i$ only and then ``blow'' it to bound a disc of symplectic area $\frac{1}{3}$, ---
it is the desired one. Therefore the moduli space of special Bohr - Sommerfeld lagrangian cycles $\tilde {\cal M}_{SBS} (S^1, 0, {\cal O}(3))$
is organaized as follows: over generic point of  $\mathbb{P}H^0 (\mathbb{C}\mathbb{P}^1, {\cal O}(3)) \backslash \Sigma$ where $\Sigma$ is the Veronose embedding of $\mathbb{C} \mathbb{P}^1$
one has three preimages, and the ramification appears over the discriminant locus, corresponding to multiple zeros, where three leaves come to one.

In a sense the presented passage from the components of sceleton to exact lagrangian submanifolds looks like the standard reduction from $\bar \partial$ - operators to hermitian connections
in the theory of stable holomorphic vector bundles, see [7]. Indeed, since the quotient space of $\bar \partial$ operators modulo locally non compact gauge group is
topologically extremely complicated one realizes the holomorphic structures by the gauge classes of hermitian connections.

The realization of special Bohr - Sommerfeld lagrangian cycles presented here via D- monotonic lagrangian submanifolds modulo Hamiltonian isotopies makes it possible to
realize  the following ``mirror symmetry dream'': in [8] one claimed that lagrangian submanifolds should correspond to vector bundles. This conjectured duality
can be realized using the moduli space of SBS lagrangian submanifolds as follows: consider in our given algrebraic variety another Lagrangian submanifold $S_0 \subset X$. Then  for any 
point of the moduli space $([S], p) \in \tilde {\cal M}_{SBS}$  take the vector space $HF(S_0; S, \mathbb{C})$ of the Floer cohomology of the pair $S_0, S$, where $S$ is a smooth 
D - monotonic lagrangian submanifold, representing the class $[S]$. Since the Floer cohomology is stable with respect to Hamiltonian isotopies,
the vector space doesn't depend on the particular choice of $S$; moreover, since the moduli space $\tilde {\cal M}_{SBS}$ is locally generated by specified
Hamiltonian isotopies this implies that globally over $\tilde {\cal M}_{SBS}$ the vector spaces are combined into a complex vector bundle, which
we denote as ${\cal F}_{S_0}$. The smoothness of the representative $S$ is important here.

Thus we get a functor from the space of lagrangian submanifolds in $X$ to the set of complex vector bundles on $\tilde{\cal M}_{SBS}$.

On the other hand the old definition of ${\cal M}_{SBS}$ as a fair subset of the direct product ${\cal B}_S \times \mathbb{P} H^0(X, L)$ admits
a strightforward introduction of a Riemannian metric on it: to do this one has to add an appropriate orientation of our Bohr - Sommerfled lagrangian
submanifolds. Indeed, since $X$ after the choice of $\omega$ is automatically endowed with the corresponding Riemannian metric $g$ fixing
an orientation on $S \in {\cal B}_S$ we get the corresponding inner product on the tangent space at each point (and if we don't fix an orientation then we get a conformal structure).
Recall that at a point $S \in {\cal B}_S$ the tangent space is given by $C^{\infty}(S, \mathbb{R})$ modulo constants (see [2]), and
in the presence of the restricted metric $g|_S$ this space can be modelled by smooth functions normalized by  $\int_S f d \mu(g|_S) = 0$.
Then the inner product is given just by the integration $\int_S f_1 f_2 d \mu(g|_S)$. Taking the standard Fubini - Study metric on the second direct summand
we induce the metric on the direct product ${\cal B}_S \times \mathbb{P} H^0(X, L)$ and then restricting it to our subset we get
a Riemannian metric on it (of course, one must check all details of the construction to ensure that one gets a correctly defined metric).

Now if we are working with the modified definition on the moduli space $\tilde{\cal M}_{SBS}$ the situation turns to be much more delicate:
how to incorporate the discussed construction to the space of classes? Is it possible? --- the answer is not quite clear.

We continue the work on these problems.

{\bf Acknowledgments.} I would like to cordially thank Ya. Eliashberg, who greatly inspired me  in my hamble studies.

$$$$

{\bf Bibliography}

[1] Nik. Tyurin, {\it ``Special Bohr - Sommerfeld lagrangian submanifolds''}, Izvestiya: Mathematics, 2016,  80:6, pp. 274–293;

[2] A. Gorodentsev, A. Tyurin, {\it ``Abelian Lagrangian Algebraic Geometry''}, Izvestiya: Mathematics, 2001, 65:3, pp. 437–467;

[3] P. Griffits, J. Harris, {\it ``Principles of algebraic geometry''}, NY, Wiley, 1978;

[4] Nik. Tyurin, arXiv:1508.06804, arXiv:1601.05975, arXiv:1609.00633;

[5] K. Cieliebak, Ya. Eliashberg, {\it ``From Stein to Weinstein and Back - Symplectic geometry of Affine Complex Geometry''}, Colloqium Publications Vol. 59, AMS (2012);
 
[6] Ya. Eliashberg, {\it ``Weinstein manifolds revisited''}, arXiv: 1707.03442; 

[7] S. Donaldson, P. Kronheimer, {\it ``The geometry of four manifolds''}, Clarendon Press, Oxford, 1990; 

[8] A. Tyurin, {\it ``Geometric quantization and Mirror symmetry''}, arXiv:math/9902027v.1.

\end{document}